\documentclass[conference]{IEEEtran}
\IEEEoverridecommandlockouts
\usepackage{cite}
\usepackage{amsmath,amssymb,amsfonts}
\usepackage{algorithmic}
\usepackage{graphicx}
\usepackage{textcomp}
\usepackage{subcaption}
\usepackage{multirow}
\usepackage{adjustbox}
\usepackage{array}
\usepackage{xcolor}
\usepackage{float}
\usepackage[left=0.680in,right=0.633in,top=0.764in,bottom=1.049in]{geometry}
\def\BibTeX{{\rm B\kern-.05em{\sc i\kern-.025em b}\kern-.08em
    T\kern-.1667em\lower.7ex\hbox{E}\kern-.125emX}}
\begin{document}
\title{Secure and Robust Communications \\for Cislunar Space Networks}
\author{
\IEEEauthorblockN{Selen Gecgel Cetin\IEEEauthorrefmark{1}, Gunes Karabulut Kurt\IEEEauthorrefmark{2}, Ángeles Vazquez-Castro\IEEEauthorrefmark{3}}
\IEEEauthorblockA{\vspace{0.3mm}
\IEEEauthorrefmark{1}Department of Electronic and Communication Engineering, Istanbul Technical University, Istanbul, Turkey\\ \vspace{0.3mm}
\IEEEauthorrefmark{2}Poly-Grames Research Center, Department of Electrical Engineering, Polytechnique Montreal, Montreal, QC, Canada\\ \vspace{0.3mm} \IEEEauthorrefmark{3}Department of Telecommunications and Systems Engineering, Autonomous University of Barcelona, Barcelona, Spain
}}
\maketitle
\begin{abstract}
There is no doubt that the Moon has become the center of interest for commercial and international actors. Over the past decade, the number of planned long-term missions has increased dramatically. This makes the establishment of cislunar space networks (CSNs) crucial to orchestrate uninterrupted communications between the Moon and Earth. However, there are numerous challenges, unknowns, and uncertainties associated with cislunar communications that may pose various risks to lunar missions.

In this study, we aim to address these challenges for cislunar communications by proposing a machine learning-based cislunar space domain awareness (SDA) capability that enables robust and secure communications. To this end, we first propose a detailed channel model for selected cislunar scenarios. Secondly, we propose two types of interference that could model anomalies that occur in cislunar space and are so far known only to a limited extent. Finally, we discuss our cislunar SDA to work in conjunction with the spacecraft communication system. Our proposed cislunar SDA, involving heuristic learning capabilities with machine learning algorithms, detects interference models with over $96\%$ accuracy. The results demonstrate the promising performance of our cislunar SDA approach for secure and robust cislunar communication.
\end{abstract}
\begin{IEEEkeywords}
Cislunar space, interference, lunar communications, machine learning, space domain awareness. 
\end{IEEEkeywords}
\section{Introduction}
Space with its unexplored and undiscovered sides has always attracted human curiosity, which is also reflected in exciting space communication projects. The first obvious attempts at space communication date back to the 1960s: Early Bird (Intelsat I) and Skynet (Telstar I) \cite{book1}. The development of space networks gained significant momentum after the 2010s, when commercial providers entered the scene. In light of the recent developments, both national and commercial players are now turning their attention beyond near-Earth space; \textit{cislunar space}. The number of lunar missions has reached unprecedented levels with more than 40 planned missions involving 80 space vehicles, by at least 10 different space agencies \cite{ref1_t1}. These missions are planned with different objectives, such as scientific discovery, lunar resource utilization, or as a forestep for deep space exploration \cite{mars_nasa, moon-mars}. Although each mission has a different focus, they mainly aim to discover the potential of the Moon.

\begin{figure}[!ht]
    \centering
    \includegraphics[width=\columnwidth]{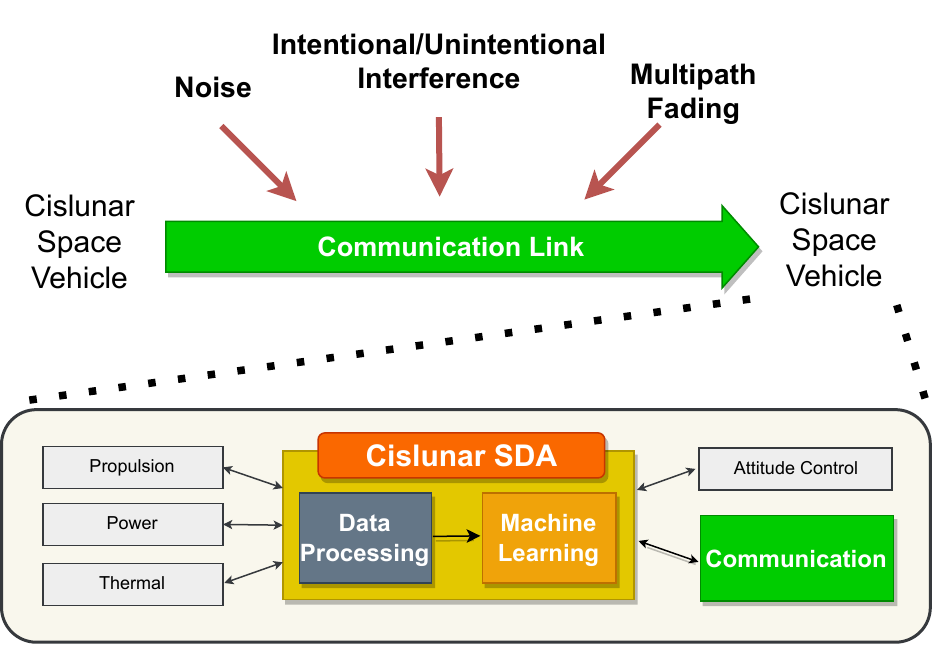}
    \caption{The diagram of the problem statement and proposed cislunar SDA system.}
    \label{fig1}
\end{figure}

When planning lunar missions, scientists and engineers must use various models to approximate real conditions, such as physical astrodynamic models, communication channel models, and atmospheric models. This is because there are many challenging issues that besiege lunar missions, such as extreme thermal environments, difficulties in communications and navigation systems, placement and maneuvering of the spacecraft, generation and storage of the power. However, uncertainties, unknowns, or unforeseen circumstances on the Moon may pose an occupational hazard for lunar missions \cite{unknowns1,unknowns2,unknowns3,unknowns4,RFI2,solarburst1,uncertain1}. The key consideration must be kept in mind that the Moon is not in close proximity to easily resolve these issues in a timely manner. Therefore, planned lunar missions must be sustainable and autonomous.

\begin{figure*}[!t]
    \centering
    \includegraphics[width=\textwidth]{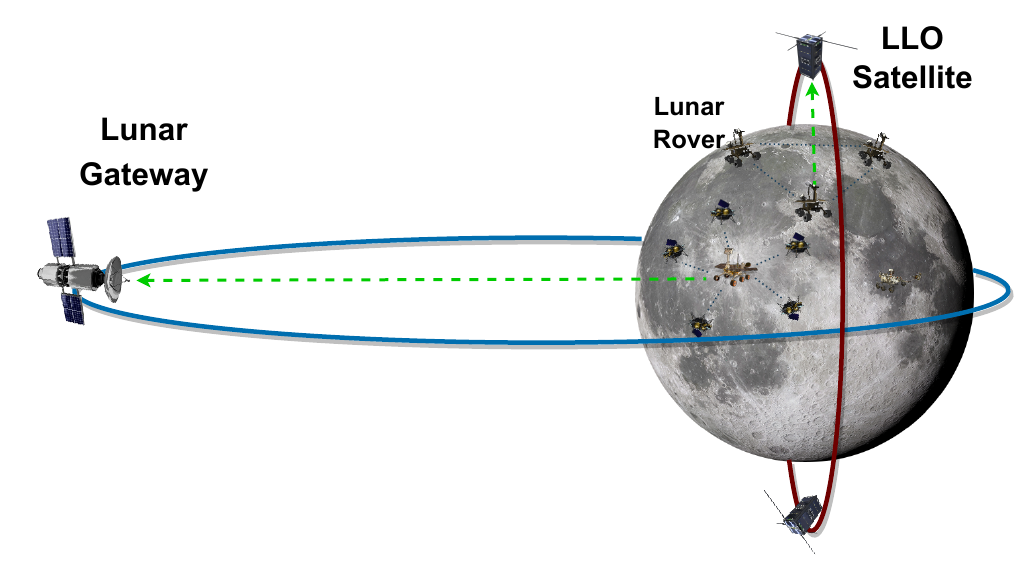}
    \caption{The overview of the two lunar communication scenarios. The first scenario considers the proximity link from a lunar surface user to the lunar gateway in near-rectilinear halo orbit (NRHO), while the second scenario considers the proximity link from a lunar surface user to a satellite in a low lunar orbit (LLO).}
    \label{fig2}
\end{figure*}
To address these constraints, future lunar missions are expected to be supported by emerging cislunar space networks (CSNs) that consist of relay satellites, orbiters and spacecrafts \cite{ref1_t1}. CSNs are expected to provide enhanced connectivity among cislunar users such as crew exploration vehicles, spacecrafts, robotic rovers, crewmember landers, and astronauts. Therefore, ensuring secure and robust communications throughout the mission plays a vital role in achieving mission objectives. While significant progress has been made on plans to establish communication networks, the reliability and security aspects of CSNs have not received adequate attention despite their importance for safe and secure operations. The vulnerability of CSNs due to wireless channel broadcast nature in remote space locations necessitates enhanced protection against intentional or unintentional threats \cite{cccds1,cccds2}.

It is evident that CSNs require advanced capabilities, including detection, tracking, and identification, to comprehend potential or current threats and be aware of their surroundings. This is known as cislunar space domain awareness (SDA), which is intended to enable better decision making with intelligent and robust applications \cite{csda1,SDAme,Falco}. As lunar networks become more crucial and the number of lunar missions increases, issues for their secure and robust communication against unknown interference characteristics emphasize the importance of cislunar SDA for CSNs. Here, machine learning algorithms must be at the core for autonomous and flexible applications of cislunar SDA, as their heuristic learning capability can adapt to unexpected situations.

The aforementioned constraints, goals, requirements, and expected capabilities cannot be realized by individual subsystems. They must be addressed in a holistic system that is as nested as possible with technology infusion to increase efficiency and reduce the cost of lunar missions. To this end, we come up with the integration of cislunar SDA and communication systems in an autonomous and intelligent manner for potential users of CSNs. Our contributions are listed below:
\begin{itemize}
    \item We develop a cislunar channel model by jointly evaluating particular channel conditions of cislunar space networks: noise, fading, and interference.
    \item We analyze the effect of thermal noise on the signal-to-noise ratio (SNR) by considering the solar insolation of the Moon during lunar phases.
    \item We propose two different interference models by considering the known and unknown circumstances in the cislunar environment.
    \item We propose an autonomous and intelligent interference detection system, as shown in Fig \ref{fig1}. Two different machine learning algorithms with high and low complexity are implemented: convolutional neural networks (CNN) and decision trees (DTREE). Both interference models are analyzed for two relevant selected scenarios of the future lunar communications architecture, as shown in Fig. \ref{fig2}. 
\end{itemize}
Overall, our proposed cislunar SDA application achieved superior detection performance over $96\%$ accuracy with both machine learning algorithms. The results also demonstrate the robustness of our solution to various behaviors of interference, such as intermittent presence with different probabilities or continuous presence, varied power levels and characteristics of the interference.

\begin{figure}[!ht]
    \centering
    \includegraphics[width=\columnwidth]{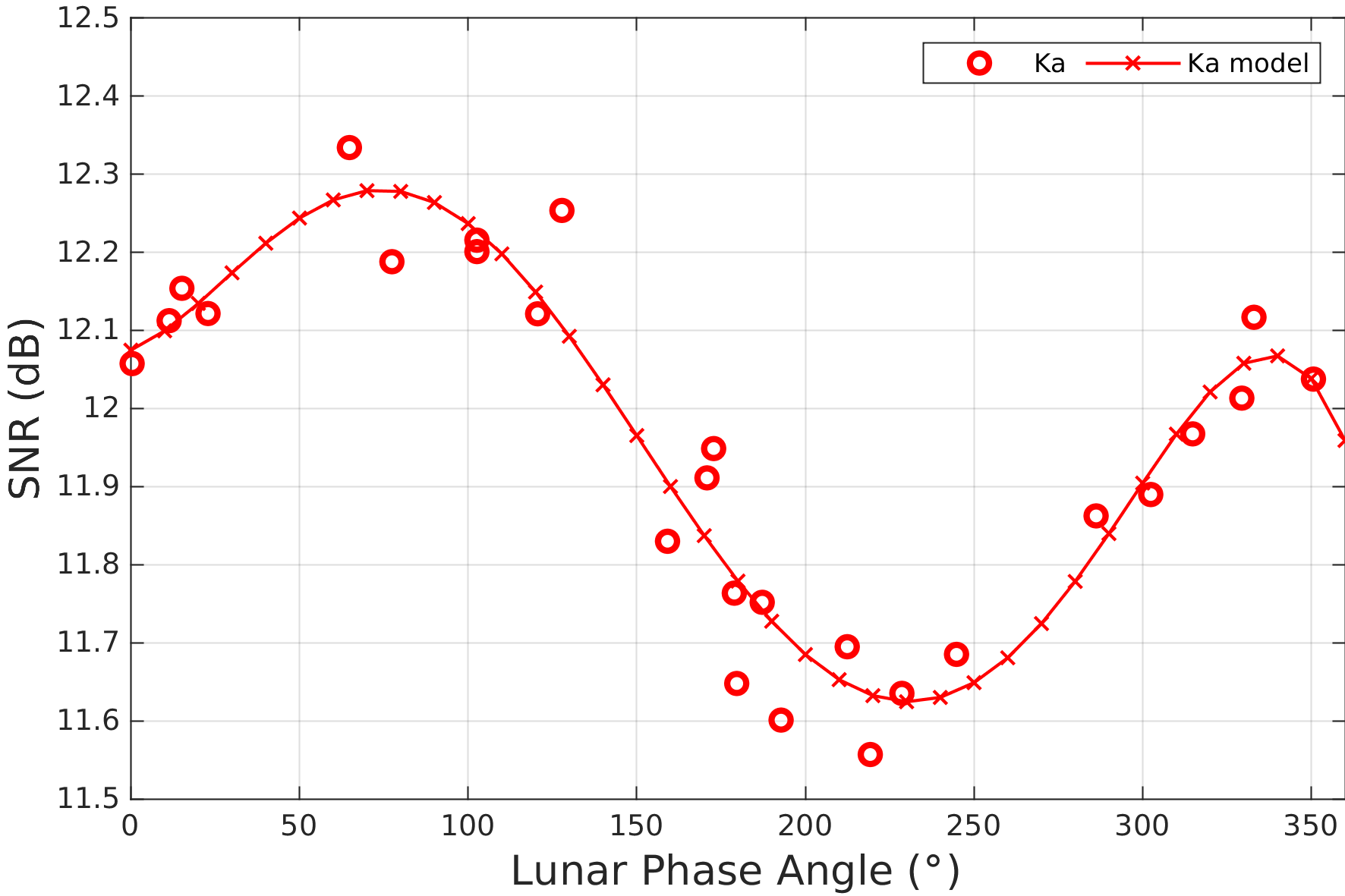}
    \caption{The impact of brightness temperature changes (depending on the lunar phase angle, $\psi$) on SNRs for the proximity link from the lunar surface user to lunar Gateway.}
    \label{Gateway_ProbStat}
\end{figure}
\begin{figure}[!ht]
    \centering
    \includegraphics[width=\columnwidth]{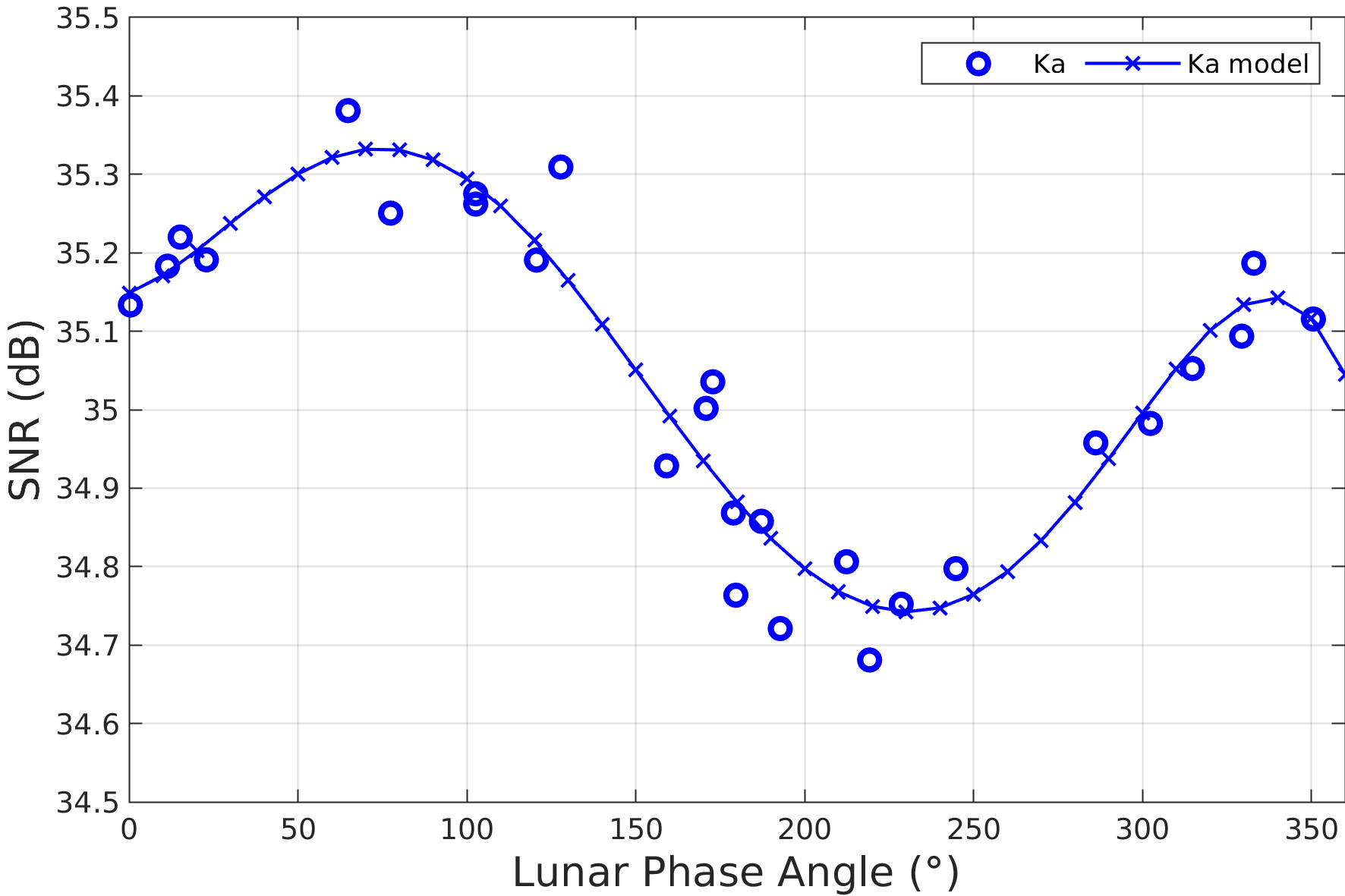}
    \caption{The impact of brightness temperature changes (depending on the lunar phase angle, $\psi$) on SNRs for the proximity link from the lunar surface user to LLO satellite.}
    \label{LLO_ProbStat}
\end{figure}

\section{Problem Statement}

In this work, we first focused on identifying the key factors affecting the cislunar communication channel: noise, multipath fading, and interference, as shown in Fig. \ref{fig1}. The noise in the communication links of CSNs may arise from electronic equipment, the cosmic microwave background and lunar noise \cite{R1}. The thermal noise is caused by thermal emissions at the receiver and generally modeled with additive white Gaussian noise (AWGN), while lunar noise is aroused from the brightness temperature of the lunar surface. The brightness temperature of the Moon shows variations depending on the lunar phases that can be defined as a function of solar illumination \cite{R7_Morabito}. The authors of \cite{R7_Morabito,R8_Morabito} provide both measurements and the approximate models of the brightness temperature depending on the lunar phase angle for different frequency bands. The variations in SNR values for Ka band are shown in Fig. \ref{Gateway_ProbStat} and Fig. \ref{LLO_ProbStat} for each scenario depending on the lunar phases when there are not any other distortion. In the remaining part of the study, the model based brightness temperature values in \cite{R7_Morabito,R8_Morabito} are used. The multipath fading is another major factor affects the signal characteristic and leads a significant variation of SNR. It must be considered for CSNs due to undulated topography of the Moon \cite{Multipath1,multipath2,RiceK2,RiceK3,RiceK4}. 

Interference is a persistent issue for space communications systems and can be unintentional (e.g., due to astronomical signals, hardware imperfections, co-channel interference) or intentional (e.g., jamming attacks, replay attacks, or spoofing attacks) \cite{Falco}. Although we have limited information about possible unintentional or intentional sources of interference due to the limited number of current missions, we know that cislunar space will soon become crowded. Given the available knowledge about past and current missions, we need to make our systems robust against interference \cite{unknowns2}. In \cite{RFI1}, the interference characteristic at the Chang’E-4 mission is presented with the single-frequency interference model and the authors have addressed the requirement for more complex interference models. The results in \cite{RFI2} show that the interference from Chang’E-3 leads the large turbulence on group delay. Although we have (not much, but limited) knowledge about these phenomena, there is no analysis yet for potential interference anomalies, multipath fading, noise, and solar insolation. Therefore, we attempt to analyze the above effects by considering their coexistence on the communication channel.
\section{Communication Link Model} %
We analyzed the proximity links between lunar surface users and cislunar orbiters in NRHO and LLO separately, as shown in Fig. \ref{fig2}. The communication link for both proximity link scenarios is modeled with the parameters in the Table \ref{Linkbudparams} as follows. The power of the received signal for our considered scenarios is defined as 
\begin{equation}
    P_r =  P_t   G_t  G_r L_{t}^{-1}L_{r}^{-1}\varphi, 
\end{equation}
where $G_t$ is the transmit antenna gain, and $G_r$ is the receiver antenna gain. $L_{t}$ and $L_{r}$ represent losses due to the transmitter and the receiver. The free space path loss is denoted by $\varphi$ and modeled as
\begin{equation}
    \varphi= \left ( \frac{4\pi fd }{c} \right )^{2},
\end{equation}
where $d$, $f$ and $c$ are the distance between the transmitter and the receiver, frequency, and speed of light, respectively. 

\begin{figure*}[!ht]
\centering
\begin{subfigure}{0.48\textwidth}
\centering
\includegraphics[width=\linewidth]{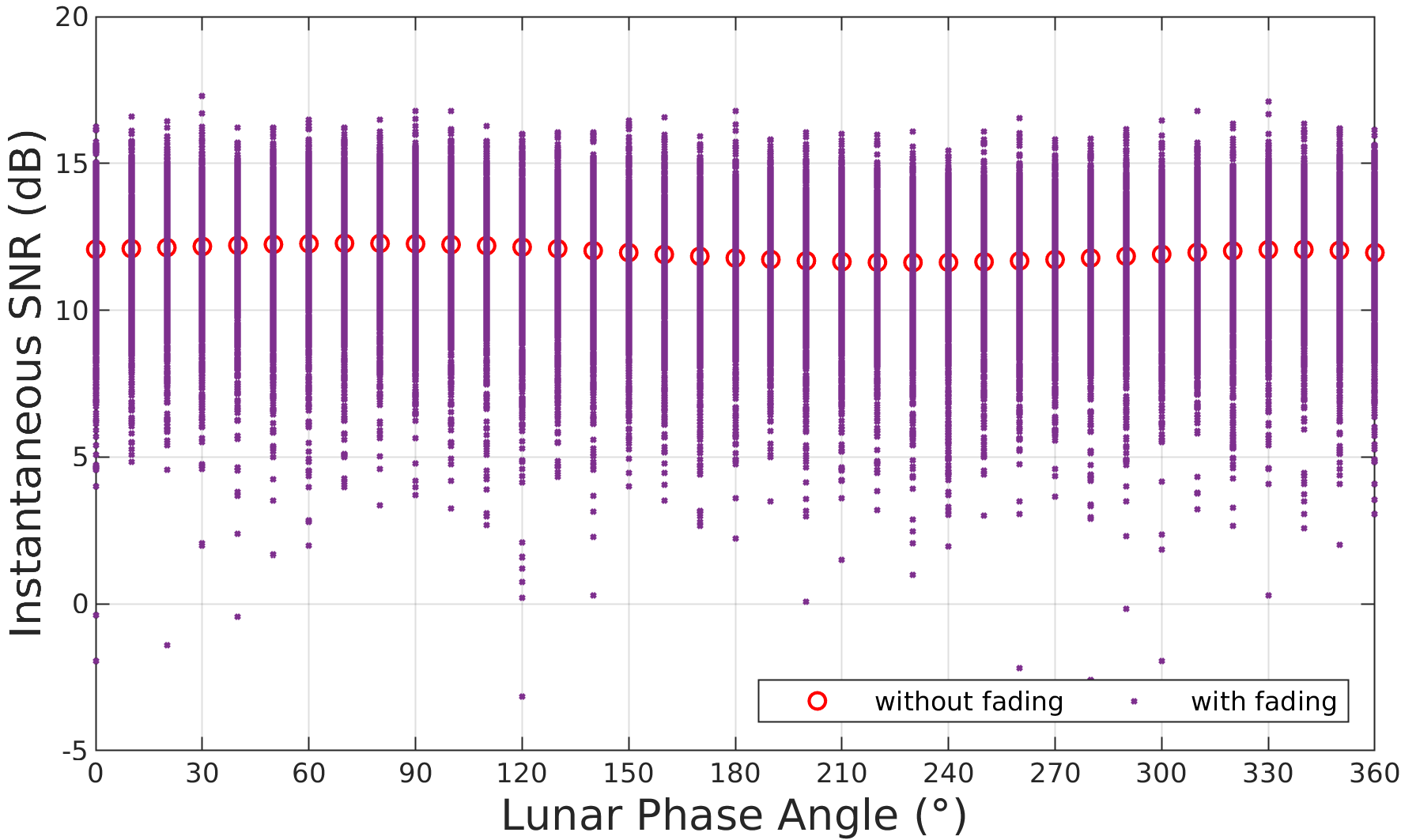}
\caption{Lunar surface user - Lunar Gateway.}
\label{Gateway_Rician}
\end{subfigure}\hspace*{\fill}
\begin{subfigure}{0.48\textwidth}
\centering
\includegraphics[width=\linewidth]{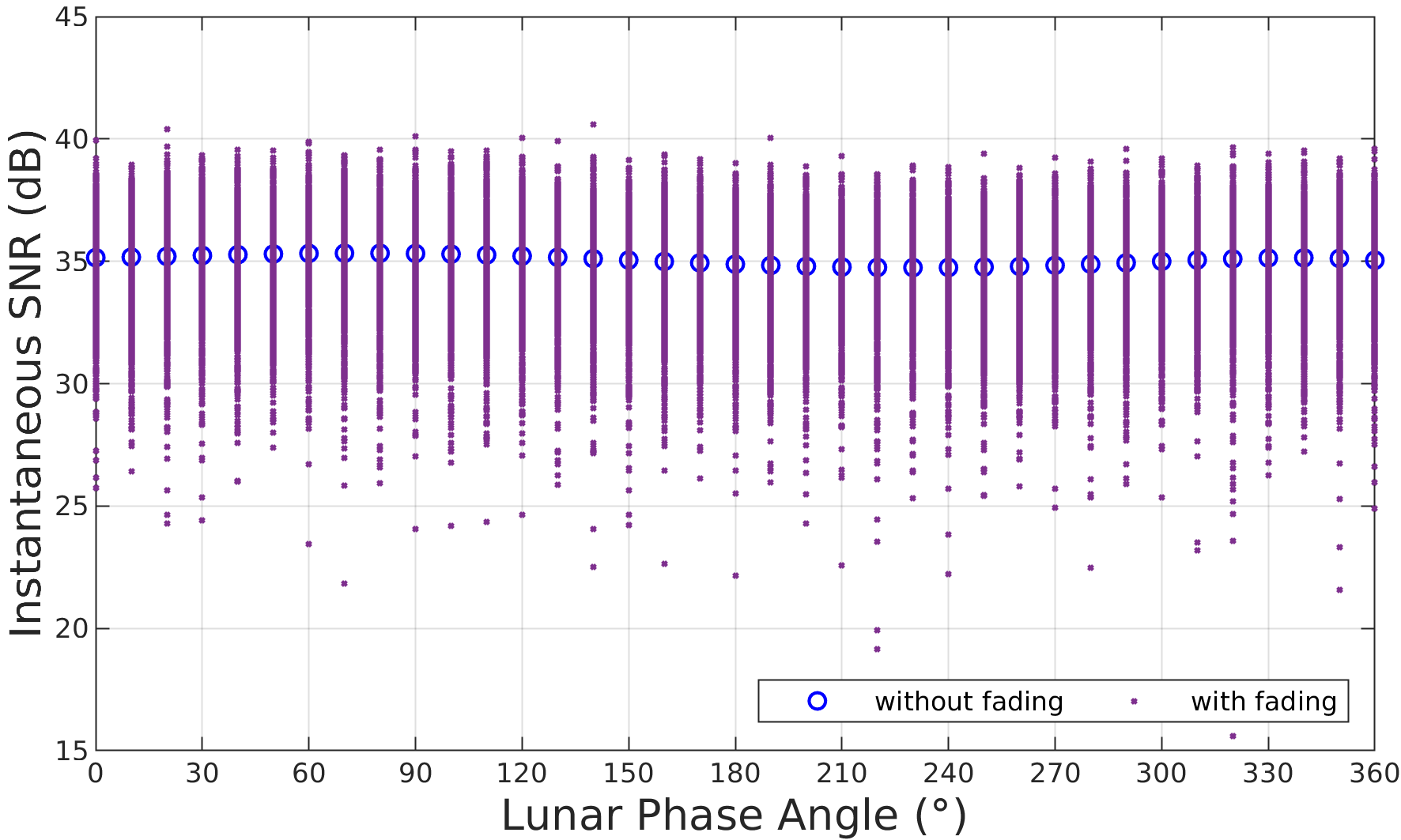}
\caption{Lunar surface user - LLO Satellite.}
\label{LLO_Rician}
\end{subfigure}
\caption{The impact of brightness temperature changes and multipath fading on the SNR for the proximity link.}
\label{int1st}
\end{figure*}
The received signal at lunar Gateway or LLO satellite is modeled by considering non-line-of-sight (NLoS) conditions \cite{Multipath1,multipath2} and AWGN as follows
\begin{equation}
    y = \sqrt{P_r}hx + n,
\end{equation}
where $h$ is channel fading coefficient with Rician distribution \cite{RiceK1,RiceK2,RiceK3,RiceK4}, and $n$ represents the additive noise follows $\mathcal{CN}(0,N_{0}(\psi))$. Here, $N_{0}(\psi)$ is one-sided power spectral density of AWGN depending on the lunar phase angle ($\psi$) and calculated as
\begin{equation}
    N_{0}(\psi) = kT_{op}(\psi)B, 
\end{equation}
where $k$ is Boltzmann’s constant $( \approx 1.38 \times 10^{-23} J/K)$, $B$ is the bandwidth in Hz, and $T_{op}(\psi)$ is the operational equivalent noise temperature in Kelvin (K). $T_{op}$ is calculated as follows 
\begin{equation}
    T_{op}(\psi) = T_{CMB} + T_{A}(\psi) + \frac{1}{\eta_{rad}}T_{TL} + \frac{1}{\eta_{rad} \eta_{TL}}T_R,
\end{equation}
where $\eta_{rad}$ and $\eta_{TL}$ show the radiation efficiency and the thermal efficiency of the transmission line. The $T_{A}(\psi)$, $ T_{TL}$, $T_R$ denote the antenna noise temperature, transmission line temperature and receiver noise temperature, respectively. $T_{CMB}$, the cosmic microwave background is approximately 2.725 K$^{{\circ}}$ \cite{R1}. The antenna noise temperature arises from physical temperature of the antenna and the brightness temperature of the subtended objects \cite{kraus}. The instantaneous received SNR is expressed as
\begin{equation}
    \gamma_{SNR} = \frac{P_r \left | h \right |^2}{N_0(\psi)},
\end{equation}
where $h$ is the channel fading random variable. The range of the instantaneous SNR are seen at different lunar phases in Fig. \ref{Gateway_Rician} and Fig \ref{LLO_Rician}.

\section{Interference Model}
We consider two different types of the interference model. The first model is characterized by considering the possibility of intermittent and continuous interference presence. Here, the intermittent presence of the interference is defined with a statistical random variable to complicate the articulations of the interference on the signal. In the second model, the interference characteristic is advanced with similar statistical behavior to noise, which makes it intractable. The proposed interference models are addressed individually through proximity links of lunar surface users with the lunar Gateway and the LLO satellite.
\begin{figure*}[!h]
\centering
\begin{subfigure}{0.48\textwidth}
\centering
\includegraphics[width=\linewidth]{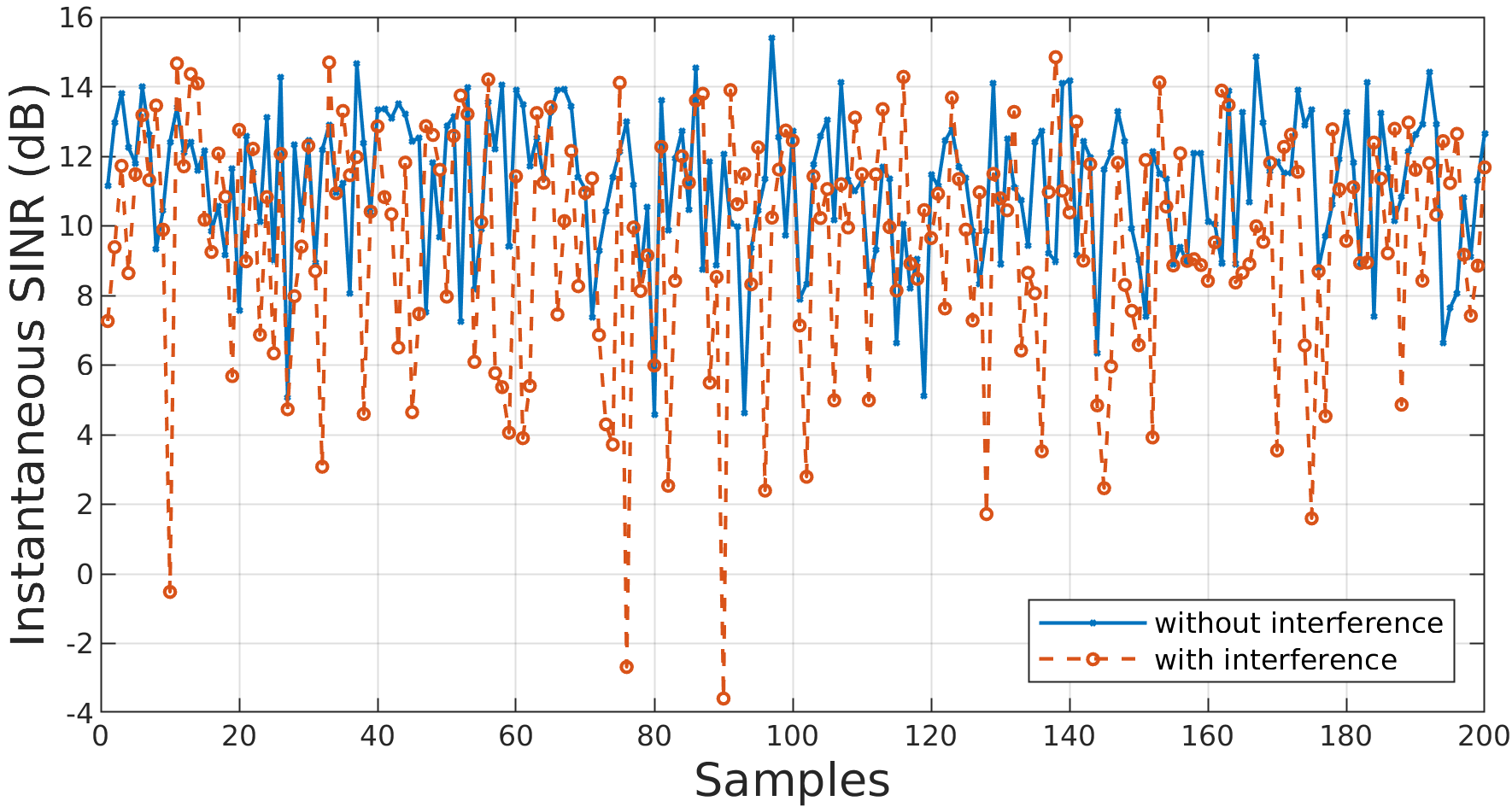}
\caption{Lunar surface user - Lunar Gateway, $\psi = 250^{\circ} $ and $p_{\alpha} = 0.5$.}
\label{Gateway_Instant_1st}
\end{subfigure}\hspace*{\fill}
\begin{subfigure}{0.48\textwidth}
\centering
\includegraphics[width=\linewidth]{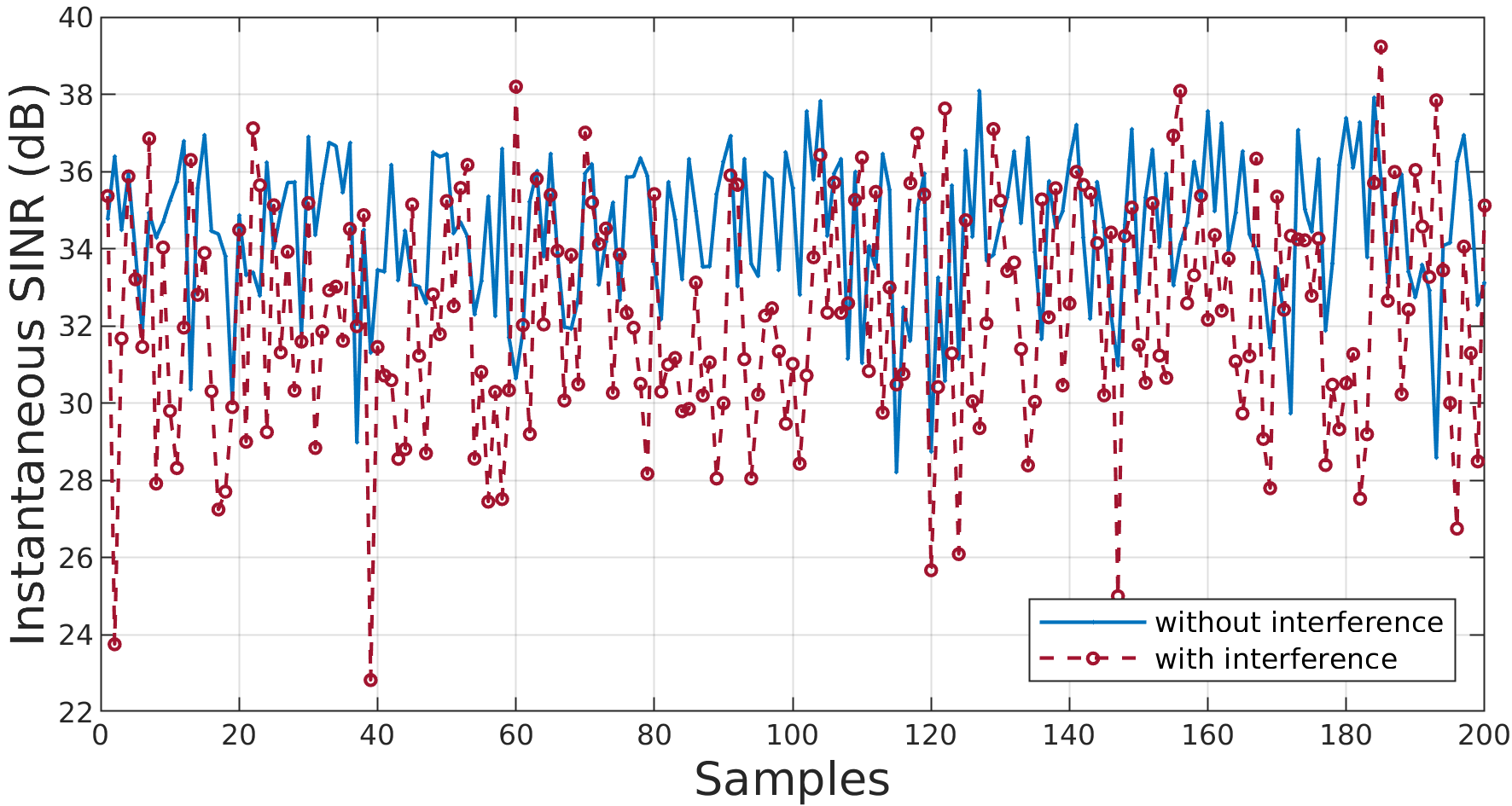}
\caption{Lunar surface user - LLO Satellite, $\psi = 240^{\circ} $ and $p_{\alpha} = 0.5$.}
\label{LLO_Instant_1st}
\end{subfigure}
\caption{The comparison of the data batches with the absence and presence of the first interference model.}
\label{int1st}
\end{figure*}
\begin{figure*}[!h]
\centering
\begin{subfigure}{0.48\textwidth}
\centering
\includegraphics[width=\linewidth]{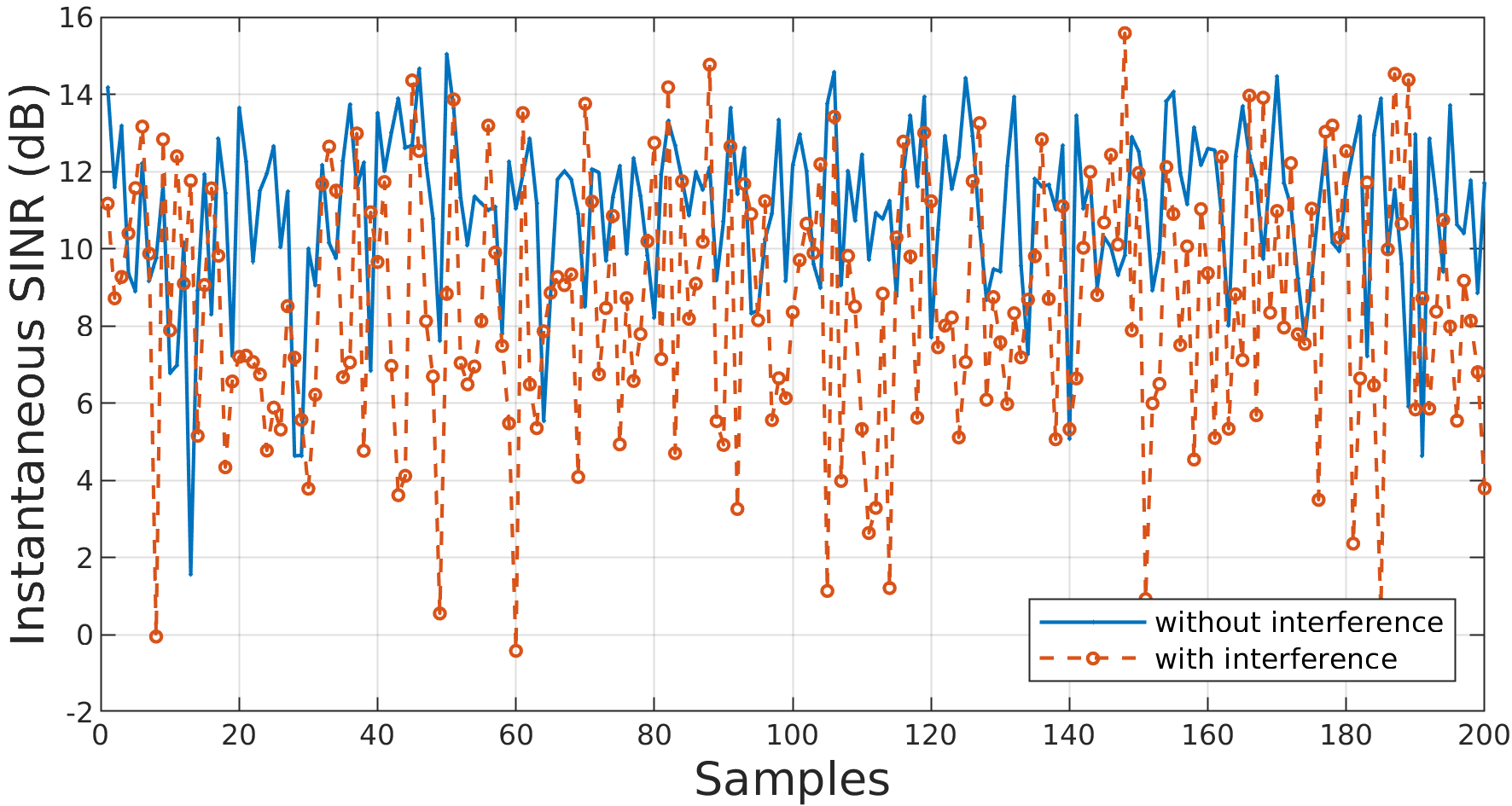}
\caption{Lunar surface user - Lunar Gateway, $\psi = 230^{\circ} $ and $p_{\alpha} = 0.5$.}
\label{Gateway_Instant_2nd}
\end{subfigure}\hspace*{\fill}
\begin{subfigure}{0.48\textwidth}
\centering
\includegraphics[width=\linewidth]{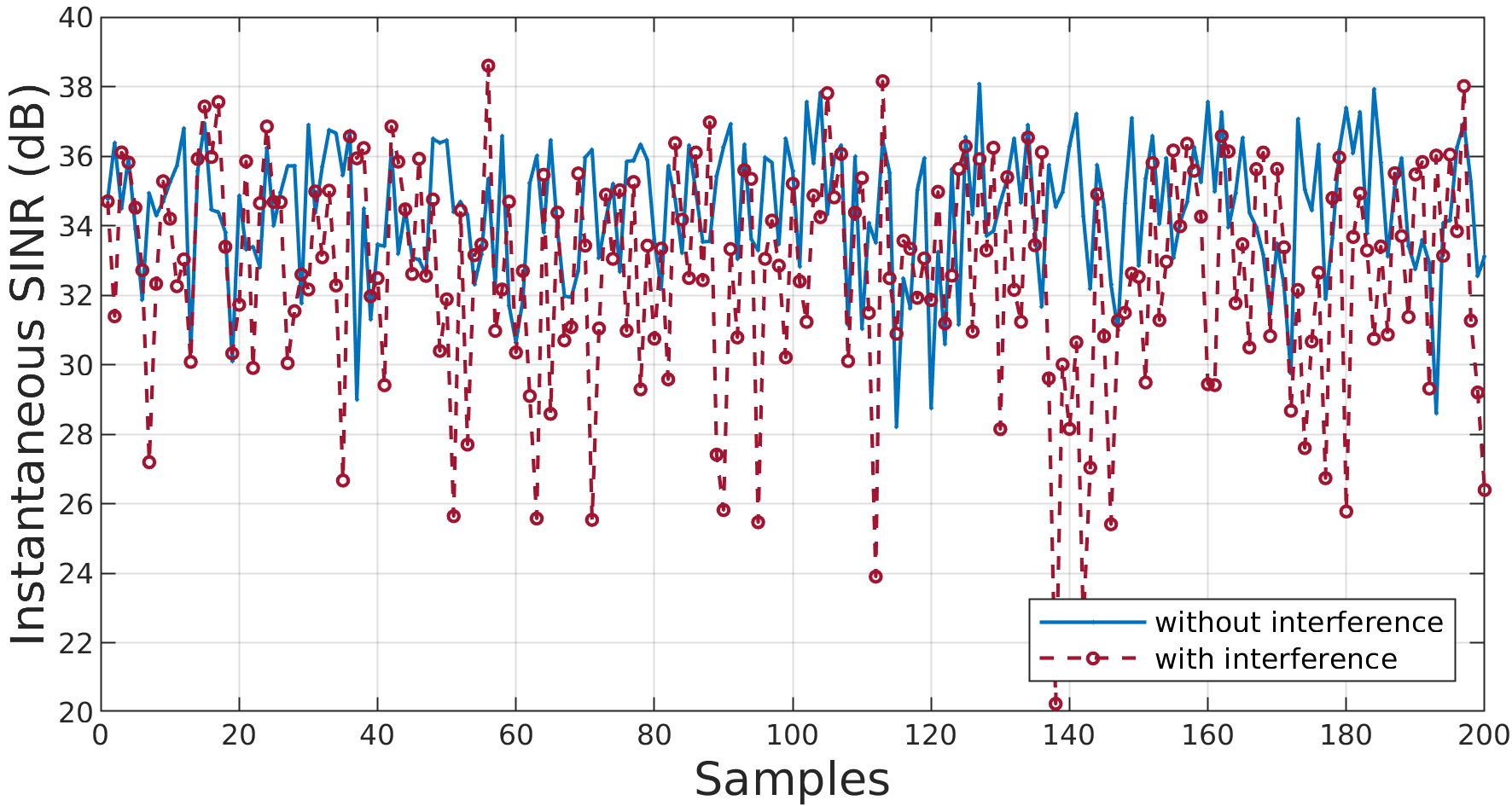}
\caption{Lunar surface user - LLO Satellite, $\psi = 240^{\circ} $ and $p_{\alpha} = 0.5$.}
\label{LLO_Instant_2nd}
\end{subfigure}
\caption{{The comparison of the data batches with the absence and presence of the second interference model.}}
\label{int2nd}
\end{figure*}
\begin{table}[!hb]
\centering
\caption{Link budget calculation parameters are used to SNR and SINR for the proposed communication links between Moon and NRHO or LLO orbiters \cite{ref1_t1,ref2_t1,ref3_t1,ref4_t1,ref5_t1,ref6_t1,ref7_t1}.}
\label{Linkbudparams}
\resizebox{\columnwidth}{!}{%
\begin{tabular}{|c|c|c|c|}
\hline
\textbf{Parameter}                                                         & \textbf{Unit} & \textbf{Moon to Gateway} & \textbf{\begin{tabular}[c]{@{}c@{}}Moon to\\ LLO Satellite\end{tabular}} \\ \hline
{Distance}                                                          & km            & 70000                    & 100                                                                      \\ \hline
Frequency                                                                  & GHz           & 26.25                    & 26.25                                                                    \\ \hline
Tx Power                                                                   & W             & 10                        & 1                                                                        \\ \hline
\begin{tabular}[c]{@{}c@{}}Tx Antenna\\ Diamater\end{tabular}              & m             & 0.254                    & 0.254                                                                    \\ \hline
\begin{tabular}[c]{@{}c@{}}Tx Antenna\\ Aperture\\ Efficiency\end{tabular} & \%            & 43                       & 43                                                                       \\ \hline
\begin{tabular}[c]{@{}c@{}}Tx Antenna\\ Gain\end{tabular}                  & dBi           & 33.21                    & 33.21                                                                         \\ \hline
Tx Losses                                                                  & dB            & 1                        & 1                                                                        \\ \hline
\begin{tabular}[c]{@{}c@{}}Rx Antenna\\ Diameter\end{tabular}              & m             & 1.5                      & 0.1                                                                      \\ \hline
\begin{tabular}[c]{@{}c@{}}Rx Antenna\\ Aperture\\ Efficiency\end{tabular} & \%            & 54                       & 56                                                                       \\ \hline
\begin{tabular}[c]{@{}c@{}}Rx Radiation\\ Efficiency\end{tabular}          & \%            & 95                       & 90                                                                       \\ \hline
\begin{tabular}[c]{@{}c@{}}Rx Antenna\\ Gain\end{tabular}                  & dBi           & 49.62                    & 26.3                                                                         \\ \hline
Rx Losses                                                                  & dB            & 3                        & 3                                                                        \\ \hline
\end{tabular}%
}
\end{table}

\subsection{Interference Model - 1}
The received signal under the first interference model is defined as
\begin{equation}
     y = \sqrt{P_r}hx + n + \sqrt{I} \alpha,
\end{equation}
where $I$ is the interference power and $\alpha \in \left \{ 0, 1 \right \}$ is status indicator of the interference follows Bernoulli distribution with the probability of $p_{\alpha}$. The instantaneous signal-to-interference-plus-noise ratio (SINR) under the first interference model is
\begin{equation}
    \gamma_{SINR} = \frac{P_r \left | h \right |^2}{N_0(\psi)+ I \left | \alpha \right |^2} . 
\end{equation}
In Fig. \ref{Gateway_Instant_1st} and Fig. \ref{LLO_Instant_1st}, the SINR of received signal under first interference model can be seen.
\subsection{Interference Model - 2}
For the second type of interference model which has similar statistical behaviour with the thermal noise, the received signal is modeled as 
\begin{equation}
     y = \sqrt{P_r}hx + n + \sqrt{I} \alpha z,
\end{equation}
where $z$ is a random variable follows $\mathcal{CN}(0,1)$. The instantaneous SINR under the second interference model is
\begin{equation}
    \gamma_{SINR} = \frac{P_r \left | h \right |^2}{ N_0(\psi) + I \left | \alpha z \right |^2 }. 
\end{equation}
In Fig. \ref{Gateway_Instant_2nd} and Fig. \ref{LLO_Instant_2nd}, the SINR of received signal under second interference model can be seen.

\section{Intelligent Detection for Cislunar SDA}
Cislunar orbiters and potential lunar surface users need research for their secure and robust transmissions, as they are essential for long-term lunar and deep space missions. In general, they are expected to be highly capable in every subsystem, from propulsion to communications. On the other hand, we know that spacecraft mass is important for missions beyond low-Earth orbit, but it is even more critical for lunar orbiters due to gravitational anomalies, placement and maneuvering difficulties. Therefore, cislunar space vehicles also be designed efficiently (in terms of fuel consumption, cost savings, and energy storage, etc.) for a long operational lifetime and in conjunction with constraints on the overall mass budget.

We propose a cislunar SDA system cooperates with subsystems of the space vehicle as shown in Fig. \ref{fig1}. In this study, we focus on the awareness of interference on the communication links against lunar Gateway and LLO satellite. We implement a low-complexity machine learning algorithm and a high-complexity one by considering the computational constraints of cislunar orbiters.
\begin{figure*}[!ht]
    \centering
    \includegraphics[width=\textwidth]{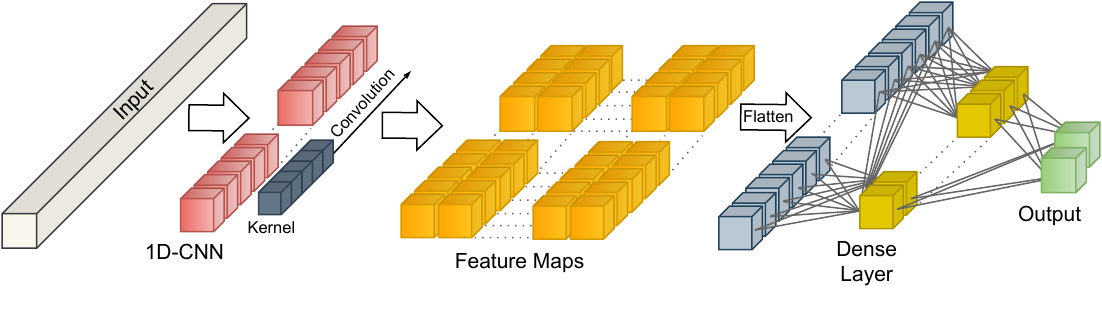}
    \caption{The proposed 1D-CNN architecture for intelligent cislunar SDA.}
    \label{fig4}
\end{figure*}
\subsection{Dataset Generation}
\label{datasetgen}
The communication signal generation is designed in accordance with the above signal models under the absence or presence of two proposed interference characteristics by using the parameters specified in the Table \ref{Linkbudparams}. The instantaneous SNR or SINR of the received signals is recorded over a time interval equal to 1000 samples and four different data sets are created. Each data set is formed to increase the detection performance and robustness of the proposed intelligent cislunar SDA system against unexpected circumstances. Therefore, all data sets consist a wide range of data under the combinations of following circumstances:
\begin{itemize}
    \item Randomly intermittent presence of (first or second) interference model with the different probabilities, $p_{\alpha} \in \left \{0.5, 0.75\right \}$.
    \item Absence of any interference model.
    \item Two particular communication scenarios in CSNs that include many differences such as; distance between transmitter and receiver, antenna features, transmit powers.
    \item The various power levels of interference models, \\ $I \in \left \{-130, -125, -120, -115, -110, -105, -100 \right \}$ in dB.
    \item Several noise levels depending on the lunar phases.
\end{itemize}
The small batches (for a time interval equal to 200 samples) of different data samples can be seen for the presence and absence of interference, each scenario and interference model in Fig. \ref{int1st} and Fig. \ref{int2nd}. Each dataset consists of 15540 training data samples and 77700 test data samples formed under above considerations.

\subsection{Machine Learning Approach}
We implemented two different machine learning algorithms with high and low computational complexity for our proposed intelligent cislunar SDA: one dimensional (1D) CNN and DTREE. CNN analyzes the patterns with the convolution operator to extract features and weights the extracted features throughout the layers. The weights of CNN models are tuned depending on the optimization algorithm and loss functions during the training process. Compared to CNN, DTREE performs estimation by constructing decision rules in a simpler strategy that aims to maximize information gain or minimize loss. In this study, DTREE model trained with 25 depths under the entropy criteria. As shown in Fig. \ref{fig4}, the CNN architecture proposed in this study consists of a one dimensional convolutional layer with 64 neurons and $5 \times 1$ kernel, a flatten layer, and 2 fully-connected layers with 128 and 2 neurons, respectively. The kernel of the CNN is randomly initialized with a uniform distribution. The rectified linear unit is used for simplification in each layer and the softmax activation function is used in the last layer. The proposed CNN architecture is optimized with the adaptive moment estimation algorithm and the categorical cross entropy loss function \cite{ADAM}.

\section{Numerical Results}
In this section, we present our results to exhibit the performance of the proposed cislunar SDA under various conditions. We firstly modeled the cislunar communication links by considering solar insolation, multipath fading, thermal noise and different types of interference. Secondly, the communication links of the lunar Gateway and the LLO satellite are analyzed individually with our models. The utilized communication links are simulated for dataset generation by including potential variations, as detailed in the Section \ref{datasetgen}. 

\begin{figure}[!hb]
\centering
\begin{subfigure}{0.495\columnwidth}
\centering
\includegraphics[width=\linewidth]{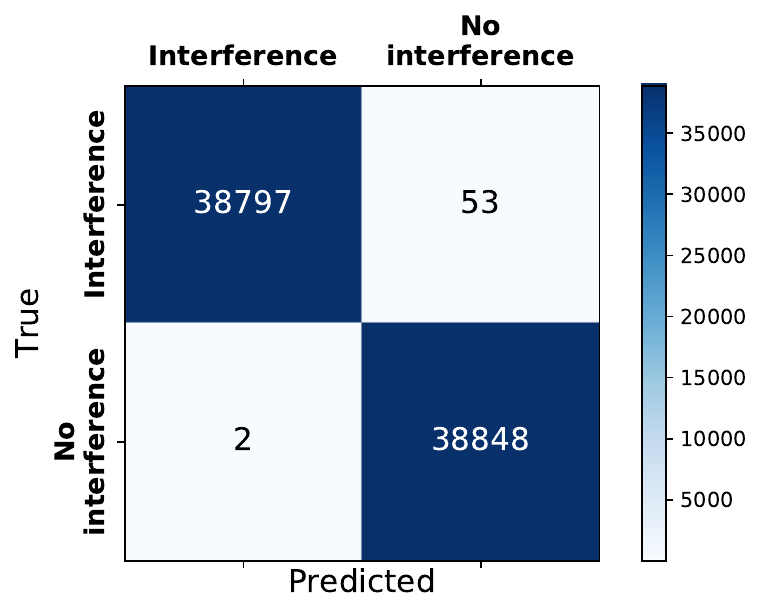}
\caption{CNN}
\label{9a1}
\end{subfigure}\hspace*{\fill}
\begin{subfigure}{0.495\columnwidth}
\centering
\includegraphics[width=\linewidth]{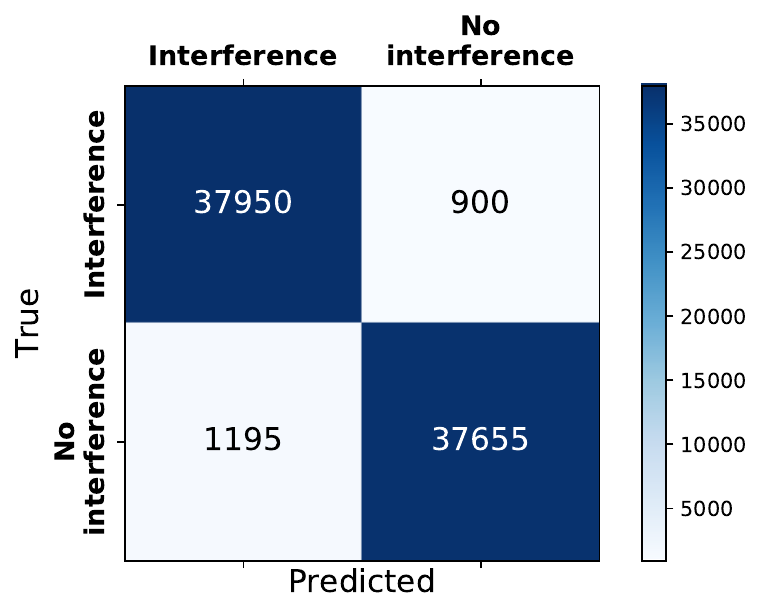}
\caption{DTREE}
\label{9b1}
\end{subfigure}\hspace*{\fill}
\caption{{Confusion matrices regarding DTREE and CNN based cislunar SDA receivers at Gateway for the first interference model.}}
\label{fig9}
\end{figure}
\begin{figure}[!ht]
\centering
\begin{subfigure}{0.495\columnwidth}
\centering
\includegraphics[width=\linewidth]{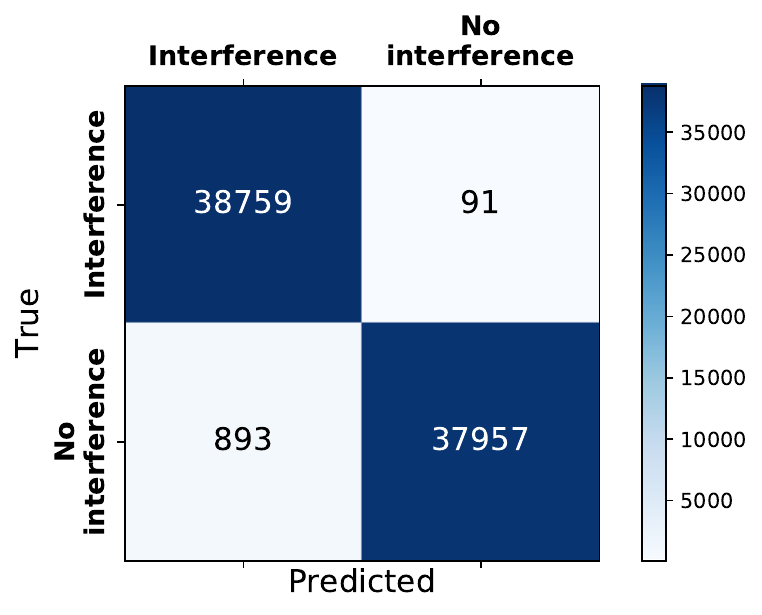}
\caption{CNN}
\label{10c1}
\end{subfigure}\hspace*{\fill}
\begin{subfigure}{0.495\columnwidth}
\centering
\includegraphics[width=\linewidth]{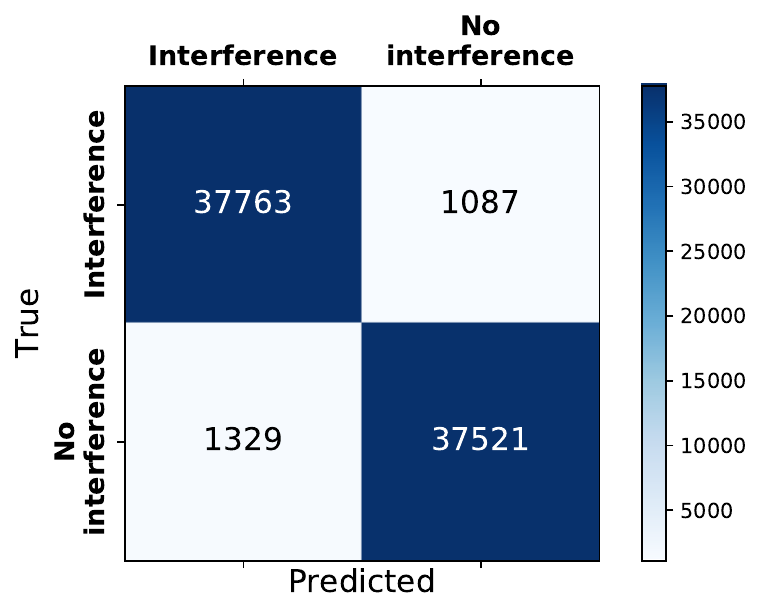}
\caption{DTREE}
\label{10d1}
\end{subfigure}
\caption{Confusion matrices regarding DTREE and CNN based cislunar SDA receivers at LLO Satellite for the first interference model.}
\label{fig10}
\end{figure}

The results in Fig. \ref{fig9} show high accuracy of $99.92\%$ and $97.30\%$ for the CNN and DTREE algorithms, respectively. The confusion matrices in Fig. \ref{9a1} and Fig. \ref{9b1} show that the precision of CNN is better than DTREE over the first interference model for the lunar Gateway link. Although the precision of CNN for the LLO satellite link decreases slightly compared to the lunar Gateway link, it shows better precision than DTREE, as shown in Fig. \ref{10c1} and Fig. \ref{10d1}. In Fig. \ref{fig10}, the detection accuracy of CNN and DTREE proves the robustness of the CNN-based cislunar SDA approach: $98.73\%$ and $96.89\%$, respectively.

\begin{figure}[!ht]
\centering
\begin{subfigure}{0.495\columnwidth}
\centering
\includegraphics[width=\linewidth]{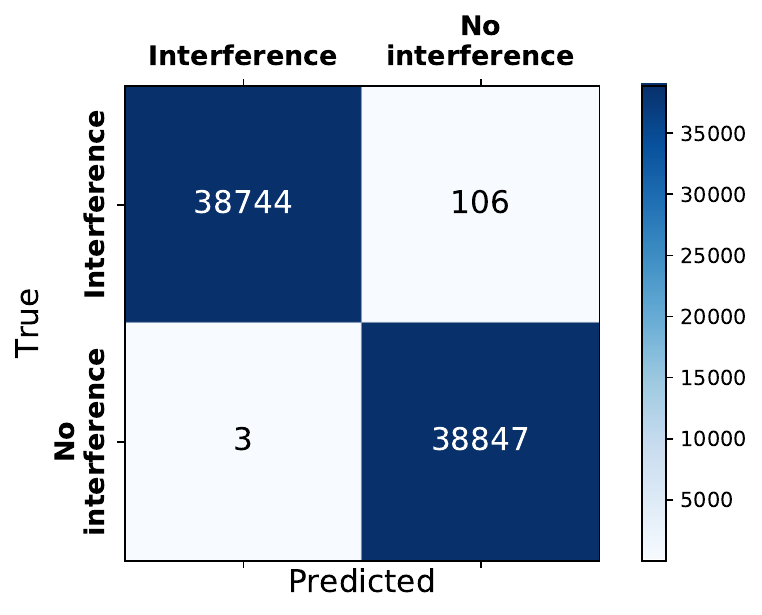}
\caption{CNN}
\label{2nd_a2}
\end{subfigure}\hspace*{\fill}
\begin{subfigure}{0.495\columnwidth}
\centering
\includegraphics[width=\linewidth]{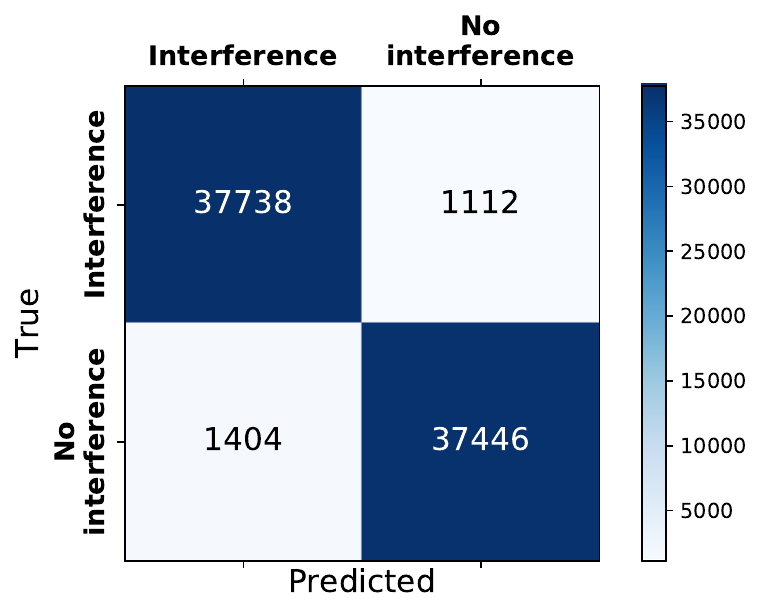}
\caption{DTREE}
\label{2nd_b2}
\end{subfigure}\hspace*{\fill}
\caption{Confusion matrices regarding DTREE and CNN based cislunar SDA receivers at Gateway for the second interference model.}
\label{fig11}
\end{figure}
\begin{figure}[!ht]
\centering
\begin{subfigure}{0.495\columnwidth}
\centering
\includegraphics[width=\linewidth]{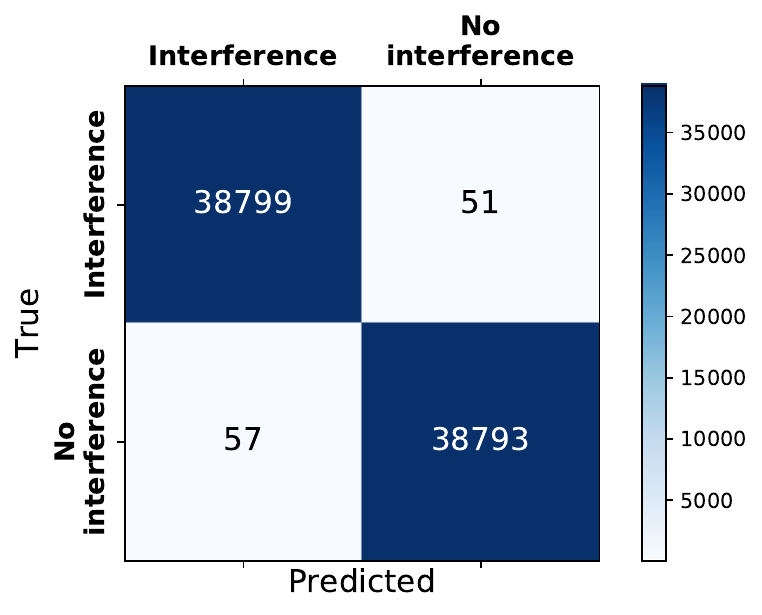}
\caption{CNN}
\label{2nd_c2}
\end{subfigure}\hspace*{\fill}
\begin{subfigure}{0.495\columnwidth}
\centering
\includegraphics[width=\linewidth]{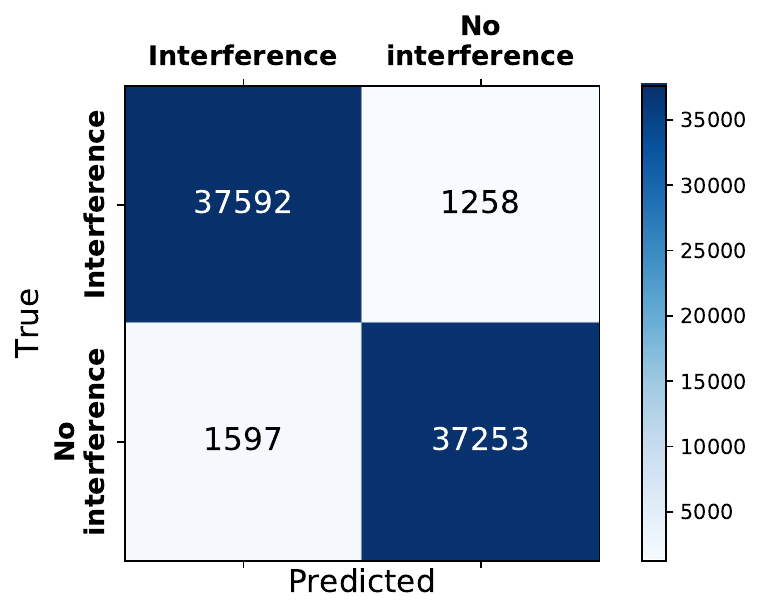}
\caption{DTREE}
\label{2nd_d2}
\end{subfigure}
\caption{Confusion matrices regarding DTREE and CNN based cislunar SDA receivers at LLO Satellite for the second interference model.}
\label{fig12}
\end{figure}

Fig. \ref{fig11} and Fig. \ref{fig12} show the superior and robust performance of the CNN over the second interference model as in the case of the first interference model. The CNN provides high detection performance with an accuracy of $99.85\%$ while the prediction accuracy of DTREE is $96.76\%$ for the detection of the first interference model and the lunar Gateway link, as shown in Fig. \ref{2nd_a2} and Fig. \ref{2nd_b2}. Although the second interference model is more challenging to detect than the first one, the precision of CNN remains stable for the lunar Gateway link, while it increases for the LLO satellite link. However, the performance of DTREE decreases slightly for both scenarios when the interference model becomes complicated. Fig. \ref{2nd_c2} and Fig. \ref{2nd_d2} show the accuracy of $99.86\%$ and $96.32\%$ for CNN and DTREE based cislunar SDA, respectively.

In summary, our proposed system with DTREE or CNN algorithm achieves a high detection accuracy of over $\%96$. As expected, we also observe that the first interference model is usually more detectable than the second interference model for the communication links of the two orbiters. In spite of the fact that our CNN-based approach shows robust performance to both interference models and its performance does not decrease. All results demonstrate the superior performance of our proposed cislunar SDA for secure and robust communication against interference and in cislunar space, and give inspiration for further but also multidisciplinary designs.
\section{Conclusions}
Lunar missions have growing goals and are becoming increasingly important. However, there are numerous challenges, unknowns, and uncertainties that pose a threat to CSNs that are designed to provide continuous communications services. As more lunar missions are conducted, the risks to CSNs become more important. For this reason, we first analyzed the corrupting effects on the communication channel and presented two interference models. Second, we proposed a cislunar SDA for secure and robust communication in cislunar space to detect interference. We have varied our solution with two machine learning algorithms and tested the performance and robustness of our approach in different cases and scenarios. The proposed cislunar SDA solution achieved superior detection performance against different interference models in several cases and scenarios. In future research, we will follow up our approach with more detailed design and anomaly models.

\bibliographystyle{IEEEtran}
\bibliography{ref}
\end{document}